\documentclass[12pt]{amsart}                   
\usepackage{amsfonts}
\usepackage{amsmath}
\usepackage{amscd}
\usepackage{amssymb}
\usepackage[dvips]{graphicx}
\usepackage[dvips]{epsfig}
\usepackage{color}

\def\?{\char'76}
\def\!{\char'74}
\def\R{\mathbb{R}}

\def\Z{\mathbb{Z}}
\def\s{\mathbb{S}^1}

\def\Im{\text{\rm Im}}

\def\max{\text{\rm max}}

\def\codim{\text{\rm codim}}

\def\coker{\text{\rm coker}}

\def\X{\widetilde{X}}

\def\B{\widetilde{B}}
\def\C{\mathcal{C}}
\def\L{\mathcal{L}}
\def\P{\mathcal{P}}
\def\U{\mathcal{U}}

\def\w{\widetilde{\omega}}
\def\prof{\text{\rm len}}

\def\strat{\mathcal{S}}

\def\m{\overline{e}}
\def\q{\overline{q}}
\def\p{\overline{p}}
\def\e{\overline{e}}
\def\top{\overline{t}}

\def\pchi{\overline{\chi}}
\def\pasito{\frac{\q-\overline{\chi}}{\q-\m}}

\def\pasitodos{\overline{2}/{1}}

\def\Omint#1#2#3{\Omega_{_{#1}}^{^{#2}}(#3)}
\def\Inv#1#2{\text{\rm I}\Omega_{_{\q}}^{^{#1}}(#2)}
\def\Homint#1#2#3{H_{_{#1}}^{^{#2}}(#3)}
\def\Omintq#1#2{\Omega_{_{\q}}^{^{#1}}(#2)}
\def\Oqchi#1#2{\Omega_{_{\q-\pchi}}^{^{#1}}(#2)}
\def\Oqm#1#2{\Omega_{_{\q-\m}}^{^{#1}}(#2)}
\def\Opasito#1#2{\Omega_{_{\pasito}}^{^{#1}}(#2)}
\def\Opasitodos#1#2{\Omega_{_{\pasitodos}}^{^{#1}}(#2)}
\def\Homintq#1#2{H_{_{\q}}^{^{#1}}(#2)}
\def\G#1#2{\mathcal{G}_{_{\q}}^{^{#1}}(#2)}
\def\GB#1{\G{#1}{B}}
\def\Cech#1#2#3{C^{^{#1}}(#2,#3)}

\def\Hcech#1#2#3{H^{^{#1}}(#2,#3)}

\def\HG#1#2{H^{^{#1}}(\mathcal{G}_{_{\q}}(#2))}
\def\HGB#1{H^{^{#1}}(\mathcal{G}_{_{\q}}(B))}
\def\Resder#1#2{\mathfrak{Low}_{_{\q}}^{^{#1}}(#2)}
\def\Resizq#1#2{\mathfrak{Upp}_{_{\q}}^{^{#1}}(#2)}
\def\Hresder#1#2{H^{^{#1}}(\mathfrak{Low}_{_{\q}}(#2))}
\def\Hresizq#1#2{H^{^{#1}}(\mathfrak{Upp}_{_{\q}}(#2))}

\def\Hpasito#1#2{\Homint{\pasito}{#1}{#2}}
\def\Hpasitodos#1#2{\Homint{\pasitodos}{#1}{#2}}
\def\Hqm#1#2{H^{^{#1}}_{_{\q-\m}}(#2)}
\def\Hqchi#1#2{H^{^{#1}}_{_{\q-\pchi}}(#2)}

\def\piarrow{\overset{\pi^{*}}\rightarrow}
\def\intarrow{\overset{\int}\rightarrow}
\def\bordarrow{\overset{\partial}\rightarrow}
\def\darrow{\overset{d}\rightarrow}

\def\eularrow{\overset{\varepsilon}\rightarrow}
\def\prarrow{\overset{pr}\rightarrow}
\def\iarrow{\overset{\imath}\rightarrow}

\def\findef{\hfill$\clubsuit$}

\newcounter{numero}
\newcommand{\Numero}{\setcounter{numero}{1}(\arabic{numero}) }
\newcommand{\numero}{\addtocounter{numero}{1}(\arabic{numero}) }

\newcounter{letra}
\newcommand{\Letra}{\medskip \setcounter{letra}{1}(\alph{letra}) }
\newcommand{\letra}{\medskip \addtocounter{letra}{1}(\alph{letra}) }

\newcounter{romnumero}

\newcounter{bibnumero}

\newtheorem{teo}{Theorem}[section]
\newtheorem{lema}[teo]{Lemma}
\newtheorem{prop}[teo]{Proposition}
\newtheorem{cor}[teo]{Corollary}

\theoremstyle{definition}
\newtheorem{definition}[teo]{}
\newtheorem{ejem}[teo]{Example}
\newtheorem{ejems}[teo]{Examples}

\theoremstyle{remark}
\newtheorem{obs}[teo]{Remark}
\newtheorem{obss}[teo]{Remarks}

\def\bteo{\begin{teo}}
\def\eteo{\end{teo}}
\def\bprop{\begin{prop}}
\def\eprop{\end{prop}}
\def\bcor{\begin{cor}}
\def\ecor{\end{cor}}
\def\blema{\begin{lema}}
\def\elema{\end{lema}}
\def\bfalsa{\begin{definition}}
\def\efalsa{\end{definition}}
\def\bejem{\begin{ejem}}
\def\eejem{\end{ejem}}
\def\bejems{\begin{ejems}}
\def\eejems{\end{ejems}}
\def\bobs{\begin{obs}}
\def\eobs{\end{obs}}
\def\bobss{\begin{obss}}
\def\eobss{\end{obss}}
\def\bdem{\begin{proof}}
\def\edem{\end{proof}}

\begin{document}

\pagestyle{myheadings} \markboth{G. Padilla}{The Gysin Sequence
for $\s$-actions on a Stratified Pseudomanifolds}

\title{The Gysin Sequence for $\s$-Actions on Stratified Pseudomanifolds}
\author{G. Padilla}
\address{Universidad Central de Venezuela-
Escuela de Matem\'atica. Caracas 1010.}
\email{gabrielp@euler.ciens.ucv.ve}
\dedicatory{ To grandmother Cira and her sisters, in loving memory.}
\date{April 24/2003}
\begin{abstract}
    For any stratified pseudomanifold $X$
    and any action of the unit circle $\s$ on $X$ preserving
    the stratification and the local structure; the orbit space $X/\s$ is also a
    stratified pseudomanifold. For each perversity $\q$ in $X$
    the orbit map $\pi:X\rightarrow X/\s$ induces a Gysin sequence
    relating the $\q$-intersection cohomologies of $X$ and $X/\s$.
    The third term of this sequence
    can be given by means of a spectral sequence on $X/\s$ whose
    second term is the cohomology of the fixed points' set $X^{\s}$
    with values on a constructible sheaf.
    The above statements
    generalize a previous work on stratified pseudomanifolds with length 1.
\end{abstract}
\maketitle

\section*{Foreword}
A stratified pseudomanifold is a topological space $X$ with two
features: the stratification and the local conical behavior. The
stratification is a decomposition of $X$ in a family of manifolds,
called strata, endowed with a partial order of incidence. The
union of open strata is a dense smooth manifold called the regular
part, its complement $\Sigma$ is the singular part of $X$. The
local conical behavior is given by the existence of charts, the
trivial model being a product $U\times c(L)$ where $U$ is a smooth
manifold and $c(L)$ is the cone of a compact stratified
pseudomanifold $L$ with lower length, we say that $L$ is a link of
$U$.\newline

When $\s$ acts on $X$ preserving the
stratification and the local structure, then the orbit space $X/\s$ is again a
stratified pseudomanifold. The orbit map $\pi:X\rightarrow X/\s$  preserves
the strata and, for each perversity $\q$ in $X$, it induces a
long exact sequence, the Gysin sequence
\[
    \cdots\rightarrow
    \Homintq{i}{X}\rightarrow
    \HG{i}{X/\s}\bordarrow
    \Homintq{i+1}{X/\s}\piarrow
    \Homintq{i+1}{X}
    \rightarrow\cdots
\]
which relates the $\q$-intersection cohomologies of $X$ and $X/\s$. The connecting homomorphism
$\partial$ depends on
the Euler class $\varepsilon\in H^{^2}_{_{\overline{2}}}(X/\s)$; it
vanishes if and only if there is a foliation on the regular part of $X$
transverse to the orbits of the action \cite{JI2}.\newline

The third complex $\G{*}{X/\s}$ in the above expression
is the Gysin term induced by the action.
Its cohomology $\HG{*}{X/\s}$ depends on basic cohomological data
of two types: local and global. There is a second long exact sequence
\[
    \cdots
    \rightarrow
    \HG{i}{X/\s}
    \rightarrow \Hresder{i}{X/\s}
    \overset{\partial'}\rightarrow
    \Hqm{i+1}{X/\s}
    \iarrow
    \HG{i+1}{X/\s}
    \rightarrow\dots
\]
where $\m$ is a perversity in $ X/\s$ vanishing on the mobile strata.
The residual term $\Hresder{*}{X/\s}$ is calculated trough a
spectral sequence whose second term
$E_{_2}=\Hcech{j}{X^{\s}}{\P^{^i}}$ is the cohomology
of the fixed points' set $X^{\s}$ with values on
a graduated constructible sheaf $\P^{^{*}}$. For each fixed point $x\in X$, the group $\s$
acts on the link $L$ of the stratum containing $x$ and the stalks{\small
\[
        \P^{^i}_{_x}=\left\{
        \begin{array}{ll}
            \Hresder{i}{L/\s} & i\leq\q(x)-3 \\
            \ker\{\partial': \Hresder{\q(x)-1}{L/\s}
        \rightarrow \Hqm{\q(x)}{L/\s}\} & i=\q(x)-2 \\
            \ker\{\partial:\HG{\q(x)-1}{L/\s}
        \rightarrow \Homintq{\q(x)+1}{L/\s}\}  & i=\q(x)-1 \\
            0 & i\geq\q(x)
        \end{array}
    \right.
\]}
are related to the Gysin sequence and the residual term of $L$.\newline


Henceforth, when we write the word {\it manifold} we are talking
about a smooth differential manifold of class
$C^{\infty}$.

\section{Stratified Spaces}
Recall the definition of a stratified spaces.
The reader will find in \cite{pflaum} a detailed exposition.
\bfalsa\label{def espacios estratificados}
    {\bf Stratified Spaces:} Let $X$ be a Hausdorff, paracompact,
    2nd countable topological space. A {\bf stratification} of $X$
    is a locally finite partition $\strat$
    whose elements are called {\bf strata}, and
    satisfy:
    \begin{itemize}
        \item[\Numero] Each stratum with the induced topology is a connected manifold.
        \item[\numero] For any two strata $S,S'\in\strat$; if $S\cap\overline{S'}$
        then $S\subset\overline{S'}$. In that case we say that $S$
    {\bf is in the border of} $S'$, and we write $S\leq S'$.
        \item[\numero] There are open strata on $X$, all of them having the
    same dimension.
    \end{itemize}
    We say that $X$ is a {\bf stratified space} whenever it has some
    stratification $\strat$.\findef
\efalsa

The open strata are also called {\bf regular}.
The {\bf regular part} of $X$ is the union of regular
strata, it is a dense open manifold.
A {\bf singular} stratum is a non-regular stratum.
The {\bf singular part} $\Sigma\subset X$ is the complement
of the regular part. \newline

The border relationship  (2) is a partial order.
A stratum is maximal (resp. minimal) if and only if it is
regular (resp. closed).
Since the stratification $\strat$ is locally finite, strict order chains
\[
    S_0<S_1<\dots<S_l
\]
on $\strat$ are always finite. The {\bf length} of $X$ is the
supreme of the integers $l$ such that there is a strict order
chain as above; we write it $\prof(X)$.

\bejems\label{ejems incidencia -ladrillo-estratrelativa}
    Here there are some examples of stratified spaces.\\
    \Numero For each manifold, the family of its connected components
    is a stratification.\\
    \numero For each stratified space $X$ and each manifold $M$, the canonical stratification
    of $M\times X$ is
    \[
        \{S\times S':S\text{ is a connected component of }M\text{ and }S'
    \text{ is a stratum of }X'\}
    \]
    \numero For each compact stratified space $L$, the
    {\bf cone of $L$} is the quotient space
    \[
        c(L)=L\times[0,\infty)/L\times\{0\}
    \]
    We write $[p,r]$ for the equivalence class of a point
    $(p,r)$, and $\star$ for the equivalence class of $L\times\{0\}$,
    which we call the {\bf vertex} of the cone. The family
    \[
        \{\star\}\sqcup\{S\times\R^{+}:S\text{ is a stratum of }L\}
    \]
    is a stratification of $c(L)$. By convention we state $c(\emptyset)=\{\star\}$.
    The {\bf radium} of the cone is the function $\rho:c(L)\rightarrow[0,\infty)$
    given  by $\rho[p,r]=r$. For each $\epsilon >0$ we write
    $c_{\epsilon}(L)=\rho^{-1}[0,\epsilon)$, it is also a stratified
    space. \findef
\eejems

Let $X$  be a stratified space. For each paracompact subspace $Y\subset X$
the {\bf induced partition} is
\[
    \strat_{Y/X}=\{C:C\text{ a connected component of }
    Y\cap S,\ S\text{ a stratum of }X\}
\]
If this family is a stratification of $Y$, then
we say that $Y$ is a {\bf stratified subspace} of $X$.\newline

A function $\alpha:X\rightarrow X'$
between two stratified spaces is a {\bf morphism}
(resp. {\bf isomorphism}) if
\begin{itemize}
    \item[\Numero] $\alpha$ is a continuous function (resp. homeomorphism).
    \item[\numero] $\alpha$ preserves the regular part, i.e.
    $\alpha(X-\Sigma)\subset(X'-\Sigma')$.
    \item[\numero] $\alpha$ sends smoothly (resp. diffeomorphically) strata into strata.
\end{itemize}
In particular, $\alpha$ is an {\bf embedding} if $\alpha(X)\subset Y$ is
open and $\alpha:X\rightarrow \alpha(X)$ is an isomorphism.
For instance, the inclusion $c_\epsilon(L)\iarrow c(L)$ is an
embedding,  and the change of radium
\[
    f:c(L)\rightarrow c_\epsilon(L)\hskip1cm
    [p,r]\mapsto [p,\epsilon\cdot\arctan(r)/\pi]
\]
is an isomorphism.

\bfalsa{\bf Stratified actions}\label{def stratified actions}
   Take a stratified space $X$, a compact abelian Lie group $G$
   and a continuous effective action
   \[
        \Phi:G\times X\rightarrow X
   \]
   We will write $\Phi(g,x)=gx$, $B=X/G$ and
   $\pi:X\rightarrow B$ for the orbit map.
   We will say that $\Phi$ is {\bf stratified} whenever
   \begin{itemize}
       \item[\Numero] The action $\Phi:G\times X\rightarrow X$ is a morphism,
       its restriction to $X-\Sigma$ is free.
       \item[\numero] For each singular stratum $S$ in $X$, the  points of $S$
       have all the same isotropy $G_S$.
   \end{itemize}
   By convention, a stratum $S$ in $X$ is {\bf mobile} (resp. {\bf fixed})
   if $G_S\neq\s$ (resp. $G_S=\s$).
\efalsa
\blema\label{lema stratified actions}
   The orbit space of a stratified action is a stratified space.
\elema
\bdem
   By \S\ref{def stratified actions}-(1), for each $g\in G$ the function $\Phi_g:X\rightarrow X$ given by
   $\Phi_g(x)=gx$ is an isomorphism. By \S\ref{def stratified actions}-(2), for each stratum $S$ of $X$
   the restriction $\pi:GS\rightarrow \pi(S)$ is a smooth locally trivial
   fibre bundle with fiber $G/G_S$.
   Since $G$ is compact;
   $B$ is a Hausdorff, paracompact, 2nd countable space \cite{bredon1}.
   The family
   \begin{equation}\label{strat}
       \strat_B=\{\pi(S):S\text{ is a stratum in }X\}
   \end{equation}
   is a locally finite partition of $B$; the border relationship passes from
   $X$ to $B$ through the orbit map. Hence $\strat_B$ is a stratification of $B$.
\edem

The family $\strat_B$ in equation (\ref{strat}) is the {\bf stratification of $B$
induced by the action $\Phi$}.
The orbit map is a morphism by construction.
\section{Stratified Pseudomanifolds and Unfoldings}
Stratified pseudomanifolds were introduced by Goresky and MacPherson in order
to extend the Poincar\'e duality to the family of stratified spaces.
The reader will find in \cite{borel} and \cite{gm1} a detailed exposition
of the subject.
\begin{definition}\label{def pseudomanifolds}
    {\bf Stratified Pseudomanifolds}
    Let $X$  be a stratified space, $S$ a stratum of $X$. A {\bf
    chart} of $S$ in $X$ is an embedding
    \[
        \alpha:U\times c(L)\rightarrow X
    \]
    where $c(L)$ is the cone of a compact stratified space, $U\subset S$ is open in $S$ and
    $\alpha(u,\star)=u $ for each $u\in U$. Notice that
    $\prof(L)<\prof(X)$. \newline

    The definition of stratified pseudomanifolds is made by
    induction on the length: we say that $X$ is
    a {\bf stratified pseudomanifold} if for each stratum $S$ there is a family  of
    charts,
    \[
        \mathcal{A}_S=
        \{\alpha:U_\alpha\times c(L)\rightarrow X\}_\alpha
    \]
    such that $\{U_\alpha\}_\alpha$ is an open cover
    of $S$, and $L$ is a compact stratified
    pseudomanifold which only depends on $S$,
    we call it a {\bf link} of $S$. \findef
\end{definition}
For any stratified pseudomanifold the link of the
regular strata is the empty set.
Any open subset of a stratified pseudomanifold is a
stratified pseudomanifold.
\bejems
    Here there are some examples of stratified pseudomanifolds.\\
    \Numero Every manifold with the canonical stratification is a
    stratified pseudomanifold.\\
    \numero For each stratified pseudomanifold $X$ and each manifold $M$, the
    product $M\times X$ is a stratified pseudomanifold.\\
    \numero If $L$ is a compact stratified pseudomanifold then $c(L)$
    is a stratified pseudomanifold.
\eejems

Recall the definition of an unfolding, we use it in order to define the
intersection cohomology of a stratified pseudomanifold
by means of differential forms.  For an introduction to the unfoldings
and their properties, see \cite{illinois}.
\begin{definition}\label{def explosiones}
    {\bf Unfoldings} Let $X$ be a stratified pseudomanifold. An
    {\bf unfolding} of $X$ is a manifold $\X$; a surjective, proper,
    continuous function
    \[
        \L:\X\rightarrow X
    \]
    and a family of unfoldings $\{\L_L:\widetilde{L}\rightarrow L\}_L$
    of the links of $X$; satisfying:
    \begin{itemize}
    \item[\Numero] The restriction $\L
    :\L^{-1}(X-\Sigma)\rightarrow X-\Sigma\ $ is a smooth trivial finite
    covering.
    \item[\numero] Each $z\in\L^{-1}(\Sigma)$ has an {\bf unfoldable chart}; i.e.,
    a commutative diagram
    \[
            \begin{array}{ccc}
                U\times \widetilde{L}\times\R &
                \overset{\tilde\alpha}{\rightarrow} &
                \X \\
                \downarrow^c & & \downarrow^{\L} \\
                U\times c(L) & \overset{\alpha}{\rightarrow} & X
            \end{array}
    \]
    where
    \begin{itemize}
        \item[\Letra] $\alpha$ is a chart.
        \item[\letra] $\widetilde\alpha$ is a diffeomorphism onto
        $\L^{-1}(\Im(\alpha))$.
        \item[\letra] Left vertical arrow is
        $c(u,\widetilde{p},t)=(u,[\L_L(\widetilde{p}),|t|])$ for each $u\in U$, $\widetilde{p}\in\widetilde{L}$,
        $t\in\R$.
    \end{itemize}
    \end{itemize}
    We say that $X$ is {\bf unfoldable} when it has an unfolding.\newline

    Let $\L:\X\rightarrow X$, $\L:\X'\rightarrow X'$ be
    two unfoldings. A morphism $\alpha:X\rightarrow X'$
    is said to be {\bf unfoldable} if there is some smooth function
    $\widetilde\alpha:\X\rightarrow \X'$ such that the following
    square
    \[
            \begin{array}{ccc}
                \X &
                \overset{\widetilde\alpha}{\rightarrow} &
                \X' \\
                \downarrow^{\L} & & \downarrow^{\L'} \\
                X & \overset{\alpha}{\rightarrow} & X'
            \end{array}
    \]
    is commutative.\findef
\end{definition}
\bfalsa{\bf Examples} These are some examples of unfoldings.\\
    \Numero For each manifold $M$ the identity $\imath:M\rightarrow M$
    is an unfolding.\\
    \numero If $\L:\X\rightarrow X$ is an unfolding then
    the product $\imath\times\L:M\times\X\rightarrow M\times X$
    is also an unfolding for any manifold $M$.\\
    \numero If $\L_L:\widetilde{L}\rightarrow L$ is an unfolding on a compact
    stratified pseudomanifold $L$; then the arrow
    $c:\widetilde{L}\times\R\rightarrow c(L)$ given by the rule $c(\widetilde{p},t)=
    [\L_L(\widetilde{p}),|t|]$ is an unfolding (cf. the left vertical arrow in
    diagram \S\ref{def explosiones}-(2)).
\efalsa

\begin{lema}\label{lema propiedades generales de las explosiones}
    Let $\L:\X\rightarrow X$ be an unfolding. Then \\
    \Numero The restriction $\L:\L^{-1}(A)\rightarrow A$ is an
    unfolding for each open subset $A\subset X$.\\
    \numero  The restriction $\L:\L^{-1}(S)\rightarrow S$
    is a smooth locally trivial fiber bundle with fiber
    $\widetilde{L}$, for each singular stratum $S$ with link $L$.
\end{lema}

\section{Intersection Cohomology}\label{seccion cohomologia de interseccion}

Between the various ways for defining the intersection
cohomology; the reader can see \cite{gm1} for a definition in
$pl$-stratified pseudomanifolds; \cite{gm2}, \cite{pervsheaves},
for a definition with sheaves; \cite{nagase} for an approach with
$\L^2$-cohomology; \cite{brylinsky} for an exposition in
Thom-Mather spaces.
In this article, we use the DeRham-like definition of intersection
cohomology as it is exposed in \cite{illinois}.
\begin{definition}\label{def formas liftables}\label{levcono}
    {\bf Liftable Forms}
    Let's fix an unfolding $\L:\X\rightarrow X$.
    A form
    $\omega\in\Omega^{^{*}}(X-\Sigma)$ is {\bf liftable} if there is
    a form $\w\in\Omega^{^{*}}(\X)$ such that $\L^{*}(\omega)=
    \w$ on ${\L^{-1}(X-\Sigma)}$.
    If such an $\w$ exists then it is
    unique by density; we call it the {\bf lifting} of
    $\omega$.  If $\omega$, $\eta$ are liftable forms then
    $d\omega$ is also liftable and we get the following equalities:
    $\widetilde{d\omega}=d\w$, $\widetilde{\omega+\eta}=\w+\widetilde\eta$,
    $\widetilde{\omega\wedge\eta}=\w\wedge\widetilde\eta$. The liftable forms
    constitute a differential graded commutative algebra.
\end{definition}
\begin{definition}\label{def cohom de interseccion}
    {\bf Intersection Cohomology}
    Let $p:M\rightarrow B$ be a surjective submersion. A smooth vector field
    $\xi$ in $M$ is {\bf vertical} if it is tangent to the fibers of $p$.
    Write $i_{\xi}$ for the contraction by $\xi$.
    The {\bf perverse degree}  $\left\| \omega \right\|_B$
    of a differential form
    $\omega\in\Omega (M)$ on $B$ is the first integer $m$ such that,
    for each vertical vector fields $\xi_0,...,\xi_m$;
    \[
        i_{\xi _0}\cdot\cdot\cdot i_{\xi_m}(\omega)=0
    \]
    Since contractions are antiderivatives of degree $-1$,
    for each $\omega,\nu\in\Omega(M)$
    \begin{equation}\label{eq grado perverso}
            \|\omega +\nu\|_B \leq \max\left\{\|\omega\|_B\ ,\|\nu\|_B\right\}
            \hskip1cm
            \|\omega\wedge\nu\|_B \leq \|\omega\|_B+\|\nu\|_B
    \end{equation}
    By convention $\left\| 0\right\| _B=-\infty$.\newline

    We define the DeRham-like intersection cohomology of $X$ by
    means of liftable differential forms and an additional
    parameter which controls their behavior when approaching to
    $\Sigma$. This new parameter is a map $\q$ which sends
    each singular stratum $S$ in $X$ to an integer $\q(S)\in\Z$; we call it
    a {\bf perversity} in $X$.\newline

    For instance, for each integer
    $n\in\Z$ we denote by $\overline{n}$ the constant perversity
    assigning $n$ to any singular stratum.
    The top perversity in $X$ is defined by $\top(S)=\codim(S)-2$
    for each singular stratum $S$.\newline

    Fix a perversity $\q$. A {\bf $\q$-form} on $X$ is a liftable form
    $\omega$ on $X-\Sigma$ satisfying
    \[
        \max\{\|\omega\|_S,\| d\omega\|_S\}\leq \q(S)
        \hskip1cm
        \forall S\text{ singular stratum}
    \]
    where, with a little abuse of language, we denote by
    $\|\omega\|_S$ the perverse degree of the restriction
    $\w\mid_{\L^{-1}(S)}$ with respect to the submersion
    $\L:\L^{-1}(S)\rightarrow S$. The $\q$-forms define a
    differential subcomplex $\Omega_{_{\q}}^{^{*}}(X)$
    whose cohomology $H_{_{\q}}^{^{*}}(X)$
    is the {\bf $\q$-intersection cohomology} of $X$.
\end{definition}

    Take two unfoldings $\L:\X\rightarrow X$ and
    $\L':\X'\rightarrow X'$. An {\bf unfoldable morphism}
    is a commutative square
    \[
            \begin{array}{ccc}
                \X &
                \overset{\widetilde\alpha}{\rightarrow} &
                \X' \\
                \downarrow^{\L} & & \downarrow^{\L'} \\
                X & \overset{\alpha}{\rightarrow} & X'
            \end{array}
    \]
    where $\alpha$ is a morphism and  $\widetilde\alpha$ is a
    smooth function, it is a {\bf lifting} of $\alpha$.
    Given a perversity $\q$ in $X'$ write still $\q$ for the perversity
    induced in $X$ in the obvious way. For each singular stratum $S$ the
    restriction of $\widetilde\alpha$ to $\L^{-1}(S)$
    preserves the perverse degree;  so $\alpha$ induces
    a well defined morphism
    \[
        \alpha^{*}:\Homintq{*}{X'}\rightarrow
    \Homintq{*}{X}
    \]

\begin{definition}\label{ejemdefhominters}
    {\bf Some properties of intersection cohomology}\\
    \Letra $H_{_{\q}}^{^{*}}(X)$ does not depend on the particular
    choice of an unfolding, for any perversity $\q$. \\
    \letra If $\q>\top$ then
    $H _{_{\q}}^{^{*}}(X)=H^{^{*}}(X-\Sigma)$ is the DeRham cohomology of $X-\Sigma$.\\
    \letra If $\q<\overline{0}$ then $H_{_{\q}}^{^{*}}(X)=H^{^{*}}(X,\Sigma)$
    is the relative cohomology of the pair. \\
    \letra If $X$ is a manifold and $\overline{0}\leq\q\leq\top$ then $H_{_{\q}}^{^{*}}(X)$
    coincides with the DeRham cohomology $H^{^{*}}(X)$.\\
    \letra For any two perversities $\p,\q$; the wedge product
    of  the forms takes into account the perversities in the following
    way:
    \[
        \Homint{\p}{i}{X}\times \Homintq{j}{X}
    \overset{\wedge}\rightarrow\Homint{\p+\q}{i+j}{X}
    \]
    In particular, the   $\overline{0}$-intersection cohomology
    $H_{_{\overline{0}}}^{^{*}}(X)$ is a differential graded algebra
    and $H_{_{\q}}^{^{*}}(X)$  is an $H_{_{\overline{0}}}(X)$-module for any
    perversity $\q$. A {\bf controlled} form is a $\overline{0}$-form.
\end{definition}

\section{Modelled Actions}

Now we introduce the family of action which we will use all
along this work, we call them modelled actions.

\begin{definition}\label{def G-pseudomanifolds}
    {\bf Modelled Actions:} Let $\Phi:G\times X\rightarrow X$
    be a stratified action.
    We say that $\Phi$ is {\bf modelled} if it satisfies
    conditions {\bf MAI} and {\bf MAII}
    stated below.\newline

    {\bf MAI} {\it For each singular stratum $S$ there is a modelled action
    $\Psi:G_S\times L\rightarrow L$ of the isotropy subgroup on the link
    of $S$.}\newline

    A {\bf modelled unfolding} of $X$ is an unfolding
    $\L:\X\rightarrow X$ in the usual sense, together with a free
    smooth action $\widetilde\Phi:G\times\X\rightarrow \X$ such that
    $\L$ is equivariant and:\\
    \Numero For each link $L$ of $X$ the induced unfolding
    $\L_L:\widetilde{L}\rightarrow L$ is modelled.\\
    \numero For each singular stratum $S$ and each
    $z\in\L^{-1}(S)$ there is a {\bf modelled chart}; i.e.,
    an unfoldable chart as in \S\ref{def explosiones}-(2), such that
    \begin{itemize}
            \item[\Letra] The diagram \S\ref{def explosiones}-(2) is
        $G_S$-equivariant.
        Here the action of $G_S$ on $U\times c(L)$ is given by
        the rule $g(u,[p,r])=(u,[gp,r])$. The free action of $G_S$
        on $U\times\widetilde{L}\times\R$ is defined as well.
        \item[\letra] The transformations of $G$ are cone-preserving:
        For each $u\in U$, $g\in G$; if
        $\Phi_g(\alpha(\{u\}\times c(L)))\cap\Im(\alpha)\neq\emptyset$
        then the arrow
            \[
                \alpha^{-1}\Phi_g\alpha\mid_u:\{u\}\times c(L)\rightarrow
                \{gu\}\times c(L)
            \]
            is an (unfoldable) isomorphism and preserves the radium
            $\rho:U\times c(L)\rightarrow[0,\infty)$.
    \end{itemize}

    {\bf MAII} {\it $X$ has a modelled unfolding.} \findef
\end{definition}

\bobs
    If $\prof(X)=0$ then $\Sigma=\emptyset$ and the conditions
    {\bf MAI}, {\bf MAII} are trivial.
    If $\prof(X)=1$ then {\bf MAI} is trivial again; condition {\bf MAII}
    can be simplified taking into account the existence of an
    equivariant normalization of the action;
    see \cite{normalizer} and \cite{gysin1}.
\eobs
\bobs
    Modelled actions constitue a category of actions. There are modelled morphisms,
    and they preserve the Euler class of the orbit space. We will have to wait untill
    \S\ref{def morfismos modelados} before we can describe them in a precise way.
\eobs

\bejems Here there are some examples of modelled actions:
    \begin{itemize}
    \item[\Numero] If $\Psi:G\times L\rightarrow L$ is a modelled action
    with a modelled unfolding $\L_L:\widetilde{L}\rightarrow L$;
    then for any manifold $U$ the induced action
    \[
        \Phi:G\times U\times c(L)\rightarrow U\times c(L)
    \hskip1cm
    g(u,[p,r])=(u,[gp,r])
    \]
    is modelled. The canonical unfolding $c:U\times\widetilde{L}\times\R
    \rightarrow U\times c(L)$ given in diagram
    \S\ref{def explosiones}-(2) is a modelled unfolding.
    \item[\numero] Let $X$ be a Thom-Mather space.
    Any stratified action $\Phi:G\times X\rightarrow X$
    preserving the tubular neighborhoods is a modelled action
    (see \cite{thom}).
    \item[\numero] If $X$ is a manifold and $\Phi$
    is a smooth effective action; then $X$ can be
    endowed with the decomposition in orbit types. This decomposition
    is a stratification and $X$ inherits an equivariant
    Thom-Mather structure. By example (1) $\Phi$ is a modelled action.
    \end{itemize}
\eejems

The main feature of modelled actions is that they preserve the category
of unfoldable pseudomanifolds.

\begin{prop}\label{prop X es Gpsve-->B es psve}
    For each modelled action $\Phi:G\times X\rightarrow X$
    \begin{itemize}
        \item[\Numero] The orbit space $B=X/G$ is a pseudomanifold.
        \item[\numero] The induced map $\L_B:\B=\X/G\rightarrow B$
    given by  the rule $\L_B(\widetilde{\pi}(x))=\pi(\L(x))$
        is an unfolding.
    \item[\numero] The orbit map $\pi:X\rightarrow B$ is
    an unfoldable morphism.
    \end{itemize}
\end{prop}
\begin{proof} Notice that, by \S\ref{lema stratified actions},
    $B$ is already a stratified
    space and $\pi:X\rightarrow B$ is a morphism.\newline

    \Numero We verify the existence of charts in $B$.
    Proceed by induction on $l=\prof(X)$,
    for $l=0$ it's trivial.
    Take some singular stratum $S$ with link $L$.
    Applying induction, by {\bf MAI}
    the quotient $L/G_S$ is a stratified pseudomanifold.
    Fix a modelled chart on $S$
    \[
            \begin{array}{ccc}
              U\times \widetilde{L}\times\R &
          \overset{\widetilde\alpha}\rightarrow &  \X \\
          \downarrow^{c} & & \downarrow^{\L} \\
              U\times c(L) & \overset{\alpha}\rightarrow
          & X
            \end{array}
    \]
    whose existence is given  by {\bf MAII}.
    Assume that $U=WV$ where $W\subset G$ is a contractible open
    neighborhood of $1\in G$, $V$ a
    slice in $S$. Write
    $\pi_L:L\rightarrow L/G_S$ for the orbit map.
    Since $\alpha$ is $G_S$-equivariant, the function
    \[
        \beta:V\times c(L/G_S)\rightarrow B
        \hskip1cm
        \beta(y,[\pi_L(p),r])=\pi\alpha(y,[p,r])
    \]
    is well defined, we will show that it is an embedding:\\
    $\bullet$ {\tt $\beta$ is injective:} Because $V$ is a slice
    in $S$ and by condition \S\ref{def G-pseudomanifolds}-(2.b), the
    transformations of $G$ are cone-preserving.\\
    $\bullet$ {\tt $\beta$ is continuous:} Because $B$, $L/G_S$
    have the respective quotient topologies.\\
    $\bullet$ {\tt $\beta$ is open:} Let $A\subset V\times
    c(L/G_S)$ be an open subset, $z\in A$. Take a compact neighborhood
    $z\in K\subset A$. Since $\beta:K\rightarrow \beta(K)$ is a continuous bijection from
    a compact space onto a Hausdorff space, it is a homeomorphism.
    There is an open $V'\subset V$ and $\epsilon >0$ such that
    \[
        z\in A'=V'\times c_{\epsilon}(L/G_S)\subset K\subset A
    \]
    So $\beta:A'\rightarrow \beta(A')$ is a homeomorphism.
    We claim that
    \[
    \beta(A')= \pi(\alpha(WV'\times c_\epsilon(L)))
    \]
    the second set is open because $WV'\subset U$ is open in $S$ and
    the orbit map $\pi$ is open. In order to sow the above equality
    take a point
    $\pi(\alpha(wv,[p,r]))\in\pi(\alpha(WV'\times c_\epsilon(L)))$.
    Then,
    \[
        w^{-1}\alpha(wv,\star)=w^{-1}(wv)=v\in\Im(\alpha)
    \]
    By condition \S\ref{def G-pseudomanifolds}-(2.b),
    $\alpha^{-1}\Phi_w^{-1}\alpha:{wv}\times c(L)\rightarrow {v}\times c(L)$
    is an isomorphism and preserves the radium. So
    \[
        \pi(\alpha(wv,[p,r]))=\pi(\alpha(v,[p',r]))=\beta(v,[\pi_L(p'),r])
    \]
    This proves that $\beta(A')\supset\pi(\alpha(WV'\times c_\epsilon(L)))$.
    The other inclusion is straightforward.\\
    $\bullet$  {\tt $\beta$ is an embedding:} On $V\times\{\star\}$
    the restriction  $\beta:V\times\{\star\}\rightarrow\pi(V)$
    given by  $\beta(y,\star)\rightarrow\pi(y)$
    is a diffeomorphism. For each stratum $R\subset L$ there is a stratum
    $S'\subset X$ such that $\alpha(V\times R\times \R^{+})\subset S'$.
    Since $\alpha$ is $G_S$-equivariant, we get the following commutative diagram
    \[
            \begin{array}{ccc}
                V\times G_S R\times\R^{+} &
                \overset{\alpha}\rightarrow & GS' \\
                \downarrow^{1\times\pi_L} & & \downarrow^\pi \\
                V\times \pi_L(R)\times\R^{+} &
                \overset{\beta}\rightarrow & \pi(S')
            \end{array}
    \]
    where the vertical arrows are submersions. So $\beta$ is
    smooth on $ V\times \pi_L(R)\times\R^{+}$ because
    $\alpha$ is smooth on $V\times G_S R\times\R^{+}$. The same argument
    can be applied to the inverse $\beta^{-1}$. Up to a change of variable,
    we can assume that $\beta(v,\star)=v\forall v$.\newline

    \numero By {\bf MAII}, the function $\L_B$ is well defined because $\L$ is equivariant.
    If $\Sigma=\emptyset$ the proof is immediate, because $\X$
    is a smooth equivariant finite trivial covering of $X$.
    In general, if $\Sigma\neq\emptyset$ then by the above remark
    $\L_B$ satisfies \S\ref{def explosiones}-(1).\newline

    Now we will prove that the charts given in the first step of this
    proof are unfoldable; this will show \S\ref{def explosiones}-(2).
    We apply induction on $l=\prof(X)$.
    For $l=0$ there is nothing to do. Take a singular stratum $S$ with link $L$.
    By induction, the $G_S$-equivariant unfolding
    $\L_L:\widetilde{L}\rightarrow L$ induces
    an unfolding $\L_{L/G_S}:\widetilde{L}/G_S\rightarrow
    L/G_S$. For each modelled chart
    $\alpha$ as in the first step of this proof, the lifted
    $\widetilde\alpha$ satisfies a smooth-like property
    analogous to \S\ref{def G-pseudomanifolds}-(2.b).
    Hence, the arrow
    \[
        \widetilde\beta:V\times\widetilde{L}/G_S\times\R
        \rightarrow\widetilde\pi(\Im(\widetilde\alpha))
        \hskip1cm
        \widetilde\beta(y,\widetilde\pi_L(\widetilde{p}),t)=
    \widetilde\pi\widetilde\alpha(y,\widetilde{p},t)
    \]
    is well defined, injective and a smooth embedding onto an open
    subset of $\B$. Notice that $\widetilde\beta$ is the
    lifting corresponding to the
    map $\beta$ given above. In consequence
    \[
            \begin{array}{ccc}
                V\times \widetilde{L}/G_S \times\R &
                \overset{\widetilde\beta}\rightarrow & \B \\
                \downarrow^{c} & & \downarrow^{\L_B} \\
                V\times c(L/G_S) &
                \overset{\beta}\rightarrow & B
            \end{array}
    \]
    is an unfoldable chart in $B$. The details are left to the reader.\newline

    \numero It is immediate from the first two statements.
\end{proof}

\section{Invariant Forms}
Some results of this \S\ where taken of \cite{S1}, \cite{coloquio santiago};
these references deal with smooth non-free circle actions on manifolds,
but the same proofs still hold in our context. The usual case of a smooth free circle
action can be seen for instance in \cite{ghv}.\newline

From now on, we fix a pseudomanifold $X$, a perversity $\q$ in $X$,
a modelled action
$\Phi:\s\times X\rightarrow X$
and a modelled unfolding $\L:\X\rightarrow X$.
\begin{definition}\label{lema quism invariants en perversas}
    {\bf Invariant forms}
    A $\q$-form $\omega$ on $X$ is {\bf invariant}
    if for each $g\in\s$ the equation $g^{*}(\omega)=\omega$ holds.
    Since
    \[
          \begin{array}{ccc}
                \X & \overset{\widetilde{\Phi_g}}\longrightarrow & \X \\
        \downarrow^{\L} & & \downarrow^{\L} \\
                X & \overset{\Phi_g}\longrightarrow & X
      \end{array}
    \]
    is an unfoldable isomorphism,
    $g^{*}:\Omega_{_{\q}}^{^{*}}(X)\rightarrow\Omega_{_{\q}}^{^{*}}(X)$
    is an isomorphism of differential complexes.
    Invariant $\q$-forms define a differential complex, denoted
    $\Inv{*}{X}$. The inclusion
    \[
        \imath:\Inv{*}{X}\rightarrow\Omintq{*}{X}
    \]
    induces an isomorphism in cohomology.
\end{definition}

Next we will study the algebraic decomposition of an
invariant $\q$-form. This decomposition depends on the existence
of an unfoldable invariant riemannian metric and an unitary
smooth vector field tangent to the orbits in the regular part.
\bfalsa\label{propexistenmetricas}
    {\bf Definitions}
    The {\bf Fundamental vector field} on $X$
    is the smooth vector field $\C$ defined on $X-\Sigma$ by
    the rule
    \[
        \C_x=d\Phi_x(\frac{\partial}{\partial g})\mid_{g=1}
    \]
    The fundamental vector field $C$ never vanishes because $X-\Sigma$ has
    no fixed points.  The lifted action
    $\widetilde\Phi:\s\times\X\rightarrow\X$ defines a
    fundamental vector field $\widetilde{\C}$ on $\X$. Notice that
    $\L_{*}(\widetilde{\C})=\C$ on $\L^{-1}(X-\Sigma)$.\newline

    An {\bf unfoldable metric} on $X$ is an invariant
    riemannian metric $\mu$ on $X-\Sigma$ such that there
    is an invariant riemannian metric $\widetilde\mu$ on $\X$ satisfying:
    \begin{itemize}
        \item[\Numero] $\L^{*}(\mu)=\widetilde{\mu}$ in $\L^{-1}(X-\Sigma)$.
        \item[\numero] $\mu\langle\C,\C\rangle=\widetilde{\mu}\langle\widetilde{\C},
        \widetilde{\C}\rangle=1$.
        \item[\numero] For each mobile stratum $S$ and each vertical vector field $\nu$
        with respect to the submersion $\L^{-1}(S)\overset{\L}{\rightarrow}S$;
        we have $\widetilde{\mu}\langle\widetilde{\C},\nu\rangle=0$.
    \end{itemize}
    For each modelled action in $X$ there is an unfoldable metric;
    the reader will find a proof in \cite{S1}.\newline

    Given an unfoldable metric $\mu$ on $X$;
    the {\bf characteristic form} induced by $\mu$  is the invariant 1-form $\chi$
    defined by the rule $\chi(v)=\mu\langle\C,v\rangle$.
\efalsa
\begin{lema}\label{cor perverse degree chi}
    The characteristic form $\chi$ satisfies
    \[
        \|\chi\|_S=\left\{
        \begin{array}{ll}
            1 & \text{$S$ a fixed stratum} \\
            0 & \text{$S$ a mobile stratum}
        \end{array}
        \right.
    \]
\end{lema}
\bdem
    By \S\ref{propexistenmetricas}-(1),
   the characteristic form $\chi$ on $X-\Sigma$ lifts to the characteristic form
   $\widetilde{\chi}$ on $\X$. The perverse degree of $\chi$
   is immediate from \S\ref{propexistenmetricas}-(2) and (3).
\edem
    Each unfoldable metric $\mu$ in $X$ induces an algebraic decomposition
    of the invariant forms. This decomposition is important in order to give
    a suitable presentation of the elements composing the Gysin Sequence of $X$.
\begin{definition}\label{def descomposicion invariant form}\label{integracion en las fibras}
    {\bf Decomposition of an Invariant Form}
    A form $\eta$ on $X-\Sigma$ is {\bf basic} if one of the following
    equivalent statements holds:
    \begin{itemize}
        \item[\Letra] $\eta$ is invariant and $\imath_{\C}(\eta)=0$.

    \item[\letra] There is a unique differential form
        $\theta$ on $B-\Sigma=\pi(X-\Sigma)$ such that $\eta=\pi^{*}(\theta)$.
    \end{itemize}
    Notice also that if $\eta$ is basic then
    \[
             0=L_{\C}(\eta)=d\imath_{\C}(\eta)+\imath_{\C}d(\eta)=
             \imath_{\C}d(\eta)
    \]
    where $L_{\C}$ is the Lie derivative with respect to the
    fundamental vector field.\newline

    For each invariant form
    $\omega\in \Inv{*}{X-\Sigma}$ there are
    $\nu\in\Omega^{^{*}}(B-\Sigma)$ and $\theta\in\Omega^{^{*-1}}(B-\Sigma)$
    satisfying
    \[
        \omega=\pi^{*}(\nu) + \chi\wedge\pi^{*}(\theta)
    \]
    The above expression is the {\bf decomposition} of $\omega$.
    The forms $\nu$, $\theta$ are uniquely determined by the
    following equations
    \[
        \pi^{*}(\theta)=\imath_{\C}(\omega)
        \hskip1cm
        \pi^{*}(\nu)= \omega-\chi\wedge\imath_{\C}(\omega)
    \]
    When $\omega$ i a liftable form then
    \[
        \w=\widetilde{\pi}^{*}(\widetilde{\nu}) +
        \widetilde{\chi}\wedge\widetilde{\pi}^{*}(\widetilde{\theta})
    \]
    So  $\nu$,$\theta$ lift respectively to $\widetilde\nu$,
    $\widetilde\theta$.\findef
\end{definition}

The first object involved in the Gysin sequence of $X$ is the
orbit map $\pi:X\rightarrow B$. For each perversity $\q$
in $X$, $\pi$ induces a well defined morphism
in intersection cohomology
\begin{equation}\label{eq morfismo pi}
    \pi^{*}:\Homintq{*}{B}\rightarrow \Homintq{*}{X}
\end{equation}
As we shall see later, the Gysin sequence is a long
exact sequence containing this map.
Now we prove that this arrow makes sense;
it is enough to study the behavior of $\pi$ with respect to
the perverse degree of an invariant $\q$-form and its algebraic
decomposition.
\begin{lema}\label{lema perverse degree invariant form}
    Take a perversity $\q$ on $X$, write also $\q$
    for the perversity induced on
    $B$ in the obvious way. Then the arrow
    \[
        \pi^{*}:\Omega_{_{\q}}^{^{*}}(B)\rightarrow
        \Inv{*}{X}
    \]
    is well defined. What's more,
    for each invariant form $\omega=\pi^{*}(\nu) +
    \chi\wedge\pi^{*}(\theta)$ and
    each singular stratum $S$, we have
    \[
        \|\omega\|_S= \max\{\|\nu\|_{\pi(S)},
        \|\chi\|_S + \|\theta\|_{\pi(S)}\}
    \]
\end{lema}
\bdem
    See \cite{S1}.
\edem
Since $\pi^{*}$ commutes with the differential $d$; the arrow
(\ref{eq morfismo pi}) is well defined.\newline

The second object that appears in the Gysin sequence of $X$ is an intersection
cohomology class in $B$ uniquely determined by the action $\Phi$; we call it
the Euler class. Up to this point, the situation is the analogous of the
smooth case.

\bfalsa\label{def eulerclass}
    {\bf The Euler class of a modelled action}
    Take an unfoldable metric $\mu$ on $X$, $\chi$ the
    characteristic form induced by $\mu$. The differential form
    $d\chi$ is basic, so there is a unique form
    $e$ on $B-\Sigma$ such that
    \[
        d\chi=\pi^{*}(e)
    \]
    This $e$ is the {\bf Euler form}
    induced by the action $\Phi$ and the metric $\mu$.
    Since $\mu$ is unfoldable, $e$ lifts to the
    Euler form $\widetilde{e}$ on $\widetilde{B}$ induced by the metric
    $\widetilde\mu$.\newline

    The {\bf Euler class} is the intersection cohomology class
    $\varepsilon=[e]\in\Homint{\e}{2}{B}$
    of the Euler form $e$ with respect to a perversity $\e$ in $B$ called the {\bf
    Euler perversity}. This $\e$ is defined by induction on the length.
    More precisely, for each singular stratum $S$ in $X$:
    \begin{itemize}
    \item[\Numero] If $S$ is mobile then we define $\e(\pi(S))=0$.
    \item[\numero] If $S$ is fixed with link $L$ and and the Euler class
    $\varepsilon_L\in\Homint{\e_L}{2}{L/\s}$ vanishes, then we define $\e(\pi(S))=1$.
    \item[\numero] If $S$ is fixed with link $L$ and
    $\varepsilon_L\in\Homint{\e_L}{2}{L/\s}\neq0$, then
    we say that $S$ is a {\bf perverse} stratum. For any perverse stratum
    $S$ we define $\e(\pi(S))=2$.
    \end{itemize}
    These properties determine the Euler perversity $\e$ in a unique way.
    Following \cite{gysin1}, the Euler class vanishes
    if and only if there is a foliation $\mathcal{F}$
    on $X-\Sigma$ transverse to the orbits of the action. \findef
\efalsa

Next we show that the Euler class is well defined.

\bprop\label{prop existen metricas buenas}
    There is an invariant unfoldable metric $\mu$ such that the Euler form
    $e$ induced by $\mu$ belongs to  $\Omint{\e}{2}{B}$.
\eprop

\bdem
We must give an unfoldable metric $\mu$ such that the induced Euler form
$e$ satisfies
\[
    \|e\|_{\pi(S)}\leq\e(\pi(S))
\]
for any singular stratum $S$ in $X$.
By \S\ref{propexistenmetricas}-(3);
we have $\|e\|_{\pi(S)}=0$ for any mobile stratum $S$.
So we only have to verify the above inequality for any fixed stratum $S$ in $X$.\newline

Proceed by induction on the length $l=\prof(X)$. For $l=0$ the action is free, so there are
no fixed strata, and the proposition trivially holds. 
We assume the inductive hypothesis, so for any fixed stratum $S$ in $X$
with link $L$; there is a metric $\mu_L$ such that the Euler form $e_L$ belongs to
$\Omint{\e}{2}{L/\s}$. The Euler class of the link 
$\varepsilon_L\in\Homint{\e_L}{2}{L/\s}$ makes sense,
as well as the classification of $S$ in perverse or non-perverse depending on the vanishing
of $\varepsilon_L$ -see \S\ref{def eulerclass}-(3). We will show that there is an unfoldable 
metric $\mu$ such that for any fixed stratum $S$ in $X$ we have
\begin{equation}\label{eq good metric}
        \|e\|_{\pi(S)}=2 \Leftrightarrow
    S\text{ is a perverse stratum}
\end{equation}
Such a metric will be called a {\it good metric}. Notice also that, by induction,
we can assume that the metric $\mu_L$ given in the link $L$ of $S$ is a good metric.
\newline

$\bullet$ {\tt Construction of a global good metric $\mu$
from a family of local ones:}\\
We give an invariant open cover
$\U=\{X_\alpha\}_\alpha$ of $X$, and a family
$\{\mu_\alpha\}_\alpha$ of unfoldable metrics such that
each $\mu_\alpha$ is a good metric in $X_\alpha$.
\begin{itemize}
    \item[\Letra] The complement of the fixed points' set
    $X_0=X-X^{\s}$ belongs to $\U$.
     We take on $X_0$ an unfoldable metric $\mu_0$.
     \item[\letra] For each fixed stratum $S$ we take a family of
     modelled charts
    \[
        \alpha:U_\alpha\times c(L) \rightarrow X
    \]
    as in \S\ref{def G-pseudomanifolds}-(2);
    such that $\{U_\alpha\}_\alpha$ is a good cover of $S$. We put
    $X_\alpha=\Im(\alpha)$ and take
    \[
        \mu_\alpha=\alpha^{-*}(\mu_{U_\alpha} + \mu_L + dr^2)
    \]
    where $\mu_{U_\alpha}$ (resp. $\mu_L$) is a riemannian (resp. good) metric
    in $U_\alpha$ (resp. in $L$). So $\mu_\alpha$
    is a good metric in $X_\alpha$.
\end{itemize}

Fix an invariant controlled
partition of the unity $\{\rho_\alpha\}_\alpha$ subordinated to $\U$. Define
\begin{equation}\label{eq mulocal}
    \mu= \underset{\alpha}\sum\ \rho_\alpha\mu_\alpha
\end{equation}

$\bullet$ {\tt Goodness of $\mu$ on a fixed stratum $S$:}
We verify the property (\ref{eq good metric}) on $S$.
\begin{itemize}
    \item[$(\Rightarrow)$]

    Write $\chi$, $e$ (resp. $\chi_\alpha$, $e_\alpha$)
    for the characteristic form and the Euler form induced  by $\mu$ on $X$
    (resp.  by $\mu_\alpha$ on $X_\alpha$). Notice that
    \begin{equation}\label{eq eulerlocal}
        d\chi=
        \underset{\alpha}\sum\ (d\rho_\alpha)\wedge\chi_\alpha
        +
        \underset{\alpha}\sum\ \rho_\alpha d\chi_\alpha
    \end{equation}
    In the above expression, the first
    sum of the right side has perverse degree $1$
    (see \S\ref{cor perverse degree chi}).
     Recall that, by
     \S\ref{lema perverse degree invariant form},
     $\|e\|_{\pi(S)}=\|d\chi\|_S$.
    If $\|d\chi\|_S=2$ then,
    by equation (\ref{eq eulerlocal}),
    \[
        \|d\chi_\alpha\|_{S\cap X_\alpha}=
    \|e_\alpha\|_{\pi(S\cap X_\alpha)}=2
    \]
    for some $X_\alpha$ intersecting $S$.
    So $\varepsilon_L\neq 0$ because $\mu_\alpha$ is a good metric.
    \item[$(\Leftarrow)$] In the rest of this proof we use some
    local properties of intersection cohomology. In particular, we
    use the step cohomology of a product $U\times c(L/\s)$ as it is
    defined in \cite{king2}. In section \S\ref{section calculos locales}
    the reader will find more details.\newline

    Assume that $\|e\|_{\pi(S)}<2$ and take
    some $X_\alpha=\Im(\alpha)\in\U$,
    the image of a modelled chart $\alpha$ on $S$.
    Write $B_\alpha=\pi(X_\alpha)\cong U_\alpha
    \times c(L/\s)$; so that
    $\|e\mid_{B_\alpha}\|_{U_\alpha}<2$.\newline

    Consider the short exact sequence of
    step intersection cohomology
    \[
        0\rightarrow \Omint{\overline{1}}{*}{B_\alpha}
    \iarrow \Omint{\overline{2}}{*}{B_\alpha}
    \overset{pr}\rightarrow
    \Opasitodos{*}{B_\alpha}
    \rightarrow 0
    \]
    which induces the long exact sequence
    \[
        \cdots\rightarrow \Homint{\overline{1}}{2}{B_\alpha}
    \rightarrow \Homint{\overline{2}}{2}{B_\alpha}
    \overset{pr^{*}}\rightarrow
    \Hpasitodos{*}{B_\alpha}
    \darrow
    \Homint{\overline{1}}{3}{B_\alpha}\rightarrow\cdots
    \]
    The inclusion $\imath_\epsilon: L/\s\rightarrow U_\alpha\times c(L/\s)$ given by
    $p\mapsto (x_0,[p,\epsilon])$, induces the isomorphism
    \[
        \imath_\epsilon^{*}:
    \Hpasitodos{2}{U_\alpha\times c(L/\s)}
    \overset{\cong}\rightarrow
    \Homint{\overline{2}}{2}{L/\s}
    \]
    where $x_0\in U_\alpha$ and $\epsilon>0$. By the above remarks,
    $(\alpha\imath_\epsilon)^{-*}(\varepsilon_L)=
    pr^{*}[e\mid_{B_\alpha}]=0$; so $\varepsilon_L=0$.
\end{itemize}
\edem

\bfalsa\label{def morfismos modelados}
    {\bf Modelled morphisms}
    Let $\Phi:\s\times X\rightarrow X$ be a modelled action.
    A {\bf perverse point} in $X$ is a point of a perverse stratum.
    We will rite $X^{perv}$  for the set of perverse  points, which is
    the union of the perverse strata, and $X^{\s}$ as usual
    for the set of fixed points.
    Let $\s\times Y\rightarrow Y$ be any other modelled action.
    An unfoldable morphism
    \[
        \begin{CD}
        \X @>{\widetilde\alpha}>> \widetilde{Y} \\
        @V{\L} VV @VV{\L'} V \\
        X @>{\alpha}>> Y\\
        \end{CD}
    \]
    is said to be {\bf  modelled} if and only each arrow in the above
    diagram is equivariant and
    \item[\Numero] $\alpha^{-1}(Y^{\s})\subset X^{\s}$.
    \item[\numero] $\alpha^{-1}(Y^{perv})\subset X^{perv}$.\newline

    In other words, a modelled morphism is an equivariant unfoldable
    morphism which preserves the clasification
    of the strata.
    For instance, any modelled chart of a fixed stratum in $X$ is a modelled
    morphism in this new sense.
\efalsa

\bteo[Functoriality of the Euler class]\label{prop funtorialidad de euler}
    The Euler class is preserved by modelled morphisms:
    If $\alpha:X\rightarrow Y$ is a modelled morphism,
    then $\alpha^{*}(\varepsilon_Y)=\varepsilon_X$.
\eteo
\bdem
    Write $\e_X,\e_Y$ for the Euler perversities in the
    orbit spaces $B_X,B_Y$. Each modelled morphism $\alpha:X\rightarrow Y$
    induces an unfoldable morphism in the orbit spaces
    \[
        \begin{CD}
        \B_X @>{\widetilde\alpha}>> \B_Y \\
        @V{\L} VV @VV{\L'} V \\
        B_X @>{\alpha}>> B_Y\\
        \end{CD}
    \]
    which we still write $\alpha$ with a little abuse of notation.
    In order to see that the arrow
    \begin{equation}\label{eq morfismo modelado}
        \alpha^{*}:\Omint{\e_Y}{*}{B_Y}\rightarrow\Omint{\e_X}{*}{B_X}
    \end{equation}
    makes sense, we consider a third perversity $\alpha^{*}[\e_Y]$
    such that the arrow
    \[
        \alpha^{*}:\Omint{\e_Y}{*}{B_Y}\rightarrow\Omint{\alpha^{*}[\e_Y]}{*}{B_X}
    \]
    is well defined and $\alpha^{*}[\e_Y]\leq\e_X$. Then
    $\Omint{\alpha^{*}[\e_Y]}{*}{B_X} \subset\Omint{\e_X}{*}{B_Y}$ and
    the map (\ref{eq morfismo modelado}) is the composition
    with the inclusion.\newline

    The perversity $\alpha^{*}[\e_Y]$ in $B_X$ is given by the rule
    \[
        \alpha^{*}[\e_Y](\pi(S))=\e_Y(\pi(R))
    \]
    for any singular strata $S,R$ respectively in $X,Y$; such that
    $\alpha(S)\subset R$. In this situation, we only need to show that
    $\alpha^{*}[\e_Y](\pi(S))\leq\e_X(\pi(S))$ or, equivalently, that
    \[
        \e_Y(\pi(R))\leq\e_X(\pi(S))
    \]
    If $R$ is mobile then $S$ is mobile because $\alpha$ is equivariant,
    so $\e_Y(\pi(R))=\e_X(\pi(S)=0$ and the inequality holds. By the other hand,
    if $R$ is fixed then the inequality is a consequence of \S\ref{def morfismos modelados},
    since $\alpha$ preserves the clasification of the strata.
\edem

\section{The Gysin Sequence}\label{section the gysin sequence}
Take a stratified pseudomanifold $X$,
a modelled action $\Phi:\s\times X\rightarrow X$ with orbit space
$B=X/\s$, and a perversity $\q$ in $X$. The orbit map
$\pi:X\rightarrow B$ preserves the strata and the perverse degree (see
\S\ref{lema perverse degree invariant form}).
Passing to the intersection cohomology we get a map
\[
    \pi^{*}:\Homintq{*}{B}\rightarrow\Homintq{*}{X}
\]
which is a string of a long exact sequence;
the Gysin sequence of $X$ induced by the action. The third
complex in the Gysin sequence is the Gysin term; its cohomology depends
at the same time on global and local basic data.
Global data concerns the Euler class $\varepsilon\in\Homint{\overline{2}}{*}{B}$ induced
by the action $\Phi$, while local data concerns the Euler classes of the
links of the perverse strata. \newline

 For instance, if $\Sigma=\emptyset$ then
$\pi:X\rightarrow B$ is a smooth $\s$-principal fiber bundle; we
get the Gysin sequence by integrating along the fibers.
If $X$ is a manifold and $\Phi$ is a smooth non-free effective action, then
$\Sigma\neq\emptyset$ and $B$ is not a manifold anymore, but a
stratified pseudomanifold. There is a Gysin sequence relating the
DeRham cohomology of $X$ with the intersection cohomology of $B$
\cite{S1}. Something analogous happens for modelled actions on
stratified pseudomanifolds with length $\leq 1$, cf. \cite{gysin1}.\newline

\begin{definition}\label{def termino de gysin}
    {\bf The Gysin Sequence}
    Fix a modelled action $\Phi:\s\times X\rightarrow X$
    on a stratified pseudomanifold $X$,
    and perversity $\q$ in $X$.
    Write still $\q$ for the obvious perversity induced on
    $B$. The orbit map $\pi:X\rightarrow B$ induces
    a short exact sequence
    \[
        0\rightarrow\Omintq{*+1}{B}
        \piarrow
        \Inv{*+1}{X}\overset{pr}\rightarrow
        \G{*}{B}\rightarrow 0
    \]
    The complex $\G{*}{B}$ is the {\bf Gysin term}.
    There is a long exact sequence
    \begin{equation}\label{eq gysinsec1}
        \dots\rightarrow
        \Homintq{i+1}{X}\overset{pr^{*}}\rightarrow
        \HGB{i}\bordarrow \Homintq{i + 2}{B}
        \piarrow
        \Homintq{i+2}{X}
        \rightarrow\dots
    \end{equation}
    This is the {\bf Gysin sequence} of $X$. \findef
\end{definition}
When the singular part of $X$ is the empty set, then $\Phi$ is a
free smooth action and $\pi:X\rightarrow B$ is a smooth principal
fiber bundle with group $\s$; so (\ref{eq gysinsec1}) is the usual
Gysin sequence. The cohomology of the Gysin term is
$\HGB{*}=\Homint{}{*}{B}$ is the DeRham cohomology and the
connecting homomorphism is the multiplication by the Euler class
(see \cite{bott}, \cite{ghv}). When $X$ has a nonempty singular
part then, by \S\ref{ejemdefhominters}, for big perversities
$\q>\top$ the sequence (\ref{eq gysinsec1}) is the usual Gysin
sequence of $X-\Sigma$ in DeRham cohomology; and for negative
perversities $\q<\overline{0}$ it is the Gysin sequence in
relative cohomology. In general;  we could naively conjecture that
$\HGB{*}=\Homintq{*}{B}$ and the connecting morphism of the Gysin
sequence to be the multiplication by the Euler class. As we will
see, the real life is richer and more complicated.\newline

Although the Gysin term is a quotient complex, it can be written
by means of basic differential forms. The {\bf characteristic
perversity} in $B$ induced by the action is
    \[
        \overline{\chi}(\pi(S))=\|\chi\|_S=\left\{
        \begin{array}{ll}
            1 & \text{$S$ a fixed stratum} \\
            0 & \text{$S$ a mobile stratum}
        \end{array}
        \right.
\]
\begin{lema}\label{lema gysin dependiendo de B}
    For each  perversity $\overline{0}\leq\q\leq\top$
    in $X$, the Gysin term $\G{*}{B}$ is isomorphic to the following
    complex{\small
    \[
    \left\{\theta\in \Oqchi{*}{B}/
        \exists\nu\in\Omega^{^{*}}(B-\Sigma):
    \begin{array}{l}
            \Numero\ \nu\text{ is liftable. }\\
        \numero\ \max\{\|\nu\|_S, \|d\nu+e\wedge\theta\|_S\}
            \leq\q(S)\forall S\text{ perverse stratum}
    \end{array}
        \right\}
    \]}
    Under this identification, the connecting homomorphism
    is
    \[
        \partial:\HGB{i}\rightarrow
        \Homintq{i}{B}\hskip1cm
        \partial[\theta]=[d\nu + e\wedge\theta]
    \]
\end{lema}
\begin{proof}
    The restriction $\pi:X-\Sigma\rightarrow B-\Sigma$ is a
    $\s$-principal fiber bundle. Consider the morphism
    of integration along the orbits
    \[
        \oint=(-1)^{^{i-1}}\pi^{-*} i_{\C}:
    \Inv{i}{X}\rightarrow\Oqchi{i-1}{B}
    \]
    defined by
    \[
        \oint\omega=(-1)^{^{i-1}}\theta
    \hskip1cm
        \omega=\pi^{*}(\nu)+\chi\wedge\pi^{*}(\theta)
    \in \Inv{i}{X-\Sigma}
    \]
    This morphism commutes with the differential.\newline

    The Gysin term
    $\G{*}{B}=\Inv{i+1}{X}/\pi^{*}(\Omintq{i+1}{B})$
    is a quotient with differential operator
    $\overline{d}(\overline{\omega})=\overline{d\omega}$, where
    $\overline{\omega}$ is the equivalence class of a differential
    form $\omega\in \Inv{i}{X}$.
    The integration along the orbits passes well to this quotient,
    \[
        \oint\G{*}{B}\rightarrow\Oqchi{*}{B}\hskip1cm
        \overline{\omega}\mapsto \oint(\omega)
    \]
    The complex given in the statement of \S\ref{lema gysin dependiendo de B}
    is the image of the above arrow. The connecting
    homomorphism arises as usual, from the Snake's Lemma.
\end{proof}

In some cases, the cohomology of the Gysin term
is closer of our naive conjecture.
\begin{prop}\label{prop simplificaciones termino de gysin}
    If $X$ has no perverse strata then
    the Euler class belongs to $H^{^2}_{_{\overline{\chi}}}(B)$
    and, for each perversity $\overline{0}\leq\q\leq\top$,
    the Gysin sequence (\ref{eq gysinsec1})
    becomes
    \[
        \dots\rightarrow
        \Homintq{i+1}{X}\intarrow
        \Hqm{i}{B} \eularrow \Homintq{i + 2}{B}
        \piarrow\Homintq{i+2}{X}
        \rightarrow\dots
    \]
    where the connecting homomorphism
    $\varepsilon$ is the multiplication by the Euler Class.
    If additionally $X$ has no fixed strata,
    then the Euler class belongs to $H^{^2}_{_{\overline{0}}}(B)$
    and the above sequence becomes
    \[
        \dots\rightarrow
        \Homintq{i+1}{X}\intarrow
        \Homintq{i}{B}\eularrow \Homintq{i + 2}{B}
        \piarrow \Homintq{i+2}{X}
        \rightarrow\dots
    \]
\end{prop}
\begin{proof}
     By \S\ref{lema gysin dependiendo de B} and the
     definition of $\pchi$, $\m$;
     the Gysin term is an intermediate complex
     \begin{equation}\label{eq gysinintermedio}
         \Oqm{*}{B}\subset\G{*}{B}\subset\Oqchi{*}{B}
     \end{equation}
    Now $X$ has no perverse strata iff $\pchi=\m$ and the extremes
    in the above inequality are identical. For the connecting homomorphism
    we take the formula in \S\ref{lema gysin dependiendo de B}
    with $\nu=0$.
\end{proof}
\begin{cor}
    If the Euler class $\varepsilon\in \Homint{\m}{2}{B}$
    vanishes then
    \[
        \Homintq{*}{X}=\Homintq{*-1}{B}\oplus
        \Hqm{*}{B}
    \]
     for each perversity $\overline{0}\leq\q\leq\top$.
     If additionally $X$ has no fixed strata, then
     \[
         H_{_{\q}}(X)=H_{_{\q}}(B)\otimes
         H(\s)
     \]
     i.e., $X$ is a cohomological product for intersection cohomology.
\end{cor}
\begin{proof}
    If the Euler class vanishes then $X$ has no perverse strata.
\end{proof}

\begin{definition}\label{def-residuosGysin}\label{def trenza de gysin}
    {\bf Residual approximations}
    Now let's assume that $X$ has perverse strata.
    We will define the residual terms; those terms allows us to
    measure the difference  between $\HGB{*}$ and the intersection cohomology
    of $B$. \newline

    Consider the inclusions in the inequality (\ref{eq gysinintermedio}).
    These inclusions provide the following short exact sequences
    \[
        \begin{array}{c}
            0\rightarrow \Oqm{*}{B}\hookrightarrow
            \GB{*}\prarrow
            \Resder{*}{B}\rightarrow 0 \\
            0\rightarrow \GB{*}\hookrightarrow
            \Oqchi{*}{B}\prarrow
            \Resizq{*}{B}\rightarrow 0
        \end{array}
    \]
    We call $\Resder{*}{B}$ (resp. $\Resizq{*}{B}$)
    the {\bf lower residue} (resp. {\bf upper residue}).
    The induced long exact sequences
    \begin{eqnarray}
            \dots\rightarrow \Hqm{i}{B}\rightarrow
            \HGB{i}\prarrow
            \Hresder{i}{B}\overset{\partial'}\rightarrow
        \Hqm{i+1}{B}\rightarrow\dots
        \label{eq resder}\\
            \dots\rightarrow \HGB{i}\rightarrow
            \Hqchi{i}{B}\prarrow
            \Hresizq{i}{B}\overset{\partial''}\rightarrow
            \HGB{i+1}\rightarrow\dots \label{eqresizq}
    \end{eqnarray}
    are the {\bf residual approximations}. Next consider the cokernel
    \[
    0\rightarrow \Oqm{*}{B}\hookrightarrow
        \Oqchi{*}{B}\prarrow
        \Opasito{*}{B}\rightarrow 0
    \]
    its cohomology $\Homint{\pasito}{*}{B}$ is called the
    {\bf step intersection cohomology} of $B$ \cite{king2}.
    The residual approximations are related by
    the long exact sequences
    \[
        \begin{array}{c}
        \dots\rightarrow \Hqm{i}{B}\rightarrow
            \Hqchi{i}{B}\rightarrow
            \Hpasito{i}{B}\darrow
            \Hqm{i+1}{B}\rightarrow\dots \\
        \dots\rightarrow \Hresder{i}{B}
        \rightarrow
            \Hpasito{i}{B}\rightarrow
            \Hresizq{i}{B}\rightarrow
        \Hresder{i}{B}\rightarrow\dots
        \end{array}
    \]
    These sequences can be arranged in a commutative exact
    diagram; called the {\bf Gysin braid} \\
    \begin{center}{\small
        \begin{picture}(30,70)
            \put(-200,40){$\Hqm{i}{B}$}         
            \put(-80,40){$\Hqchi{i}{B}$}
            \put(20,40){$\Hresizq{i}{B}$}
            \put(140,40){$\Hresder{i+1}{B}$}
            \put(-150,0){$\HGB{i}$}             
            \put(-30,0){$\Hpasito{i}{B}$}
            \put(85,0){$\HGB{i+1}$}
            \put(-210,-40){$\Hresizq{i-1}{B}$}         
            \put(-100,-40){$\Hresder{i}{B}$}
            \put(40,-40){$ \Hqm{i+1}{B}$}
            \put(165,-40){$ \Hqchi{i+1}{B}$}
            \thinlines                               
            \qbezier(-175,55)(-125,80)(-75,55)      
            \put(-80,55){\line(1,0){5}}
            \put(-75,55){\line(0,1){5}}
            \qbezier(-50,55)(0,80)(50,55)          
            \put(45,55){\line(1,0){5}}
            \put(50,55){\line(0,1){5}}
            \qbezier(70,55)(125,80)(180,55)        
            \put(175,55){\line(1,0){5}}
            \put(180,55){\line(0,1){5}}
            \qbezier(-175,-55)(-125,-80)(-75,-55)      
            \put(-80,-55){\line(1,0){5}}
            \put(-75,-55){\line(0,-1){5}}
            \qbezier(-50,-55)(0,-80)(50,-55)          
            \put(45,-55){\line(1,0){5}}
            \put(50,-55){\line(0,-1){5}}
            \qbezier(70,-55)(125,-80)(180,-55)        
            \put(175,-55){\line(1,0){5}}
            \put(180,-55){\line(0,-1){5}}
            \put(-165,30){\vector(2,-1){30}}            
            \put(-50,30){\vector(2,-1){30}}
            \put(70,30){\vector(2,-1){30}}
            \put(-110,-10){\vector(2,-1){30}}
            \put(10,-10){\vector(2,-1){30}}
            \put(130,-10){\vector(2,-1){30}}
            \put(-165,-20){\vector(2,1){30}}            
            \put(-55,-25){\vector(2,1){30}}
            \put(70,-20){\vector(2,1){30}}
            \put(-105,17){\vector(2,1){30}}
            \put(5,17){\vector(2,1){30}}
            \put(130,20){\vector(2,1){30}}
%
%
        \end{picture}\\
        \vskip2cm
        \ \\}
    \end{center}
\end{definition}

\section{The Gysin Theorem}\label{section calculos locales}
We devote the rest of this work to calculate the residual
cohomologies. The final goal is
to relate $\HGB{*}$ with basic local cohomological data by
means of the residual approximations; so we start this \S\ with
the local properties of the residues. An introduction to
presheaves, sheaves and Cech cohomology  can be found in
\cite{bott}, \cite{godement}; some results of this \S\  were taken
from \cite{illinois}.\newline

Recall that a {\bf presheaf} $\P$ on $X$ is {\bf complete} (or it is a {\bf sheaf}) iff, for any
open cover $\U=\{X_\alpha\}_\alpha$ of $X$, the augmented Cech differential
complex
\begin{equation}\label{eq resolucion de cech}
        0\rightarrow
        \P(X)\overset{\delta}{\rightarrow}
    C^{^0}(\U,\P)
    \overset{\delta}\rightarrow
        C^{^1}(\U,\P)\overset{\delta}\rightarrow
        C^{^2}(\U,\P)
    \overset{\delta}\rightarrow\cdots
\end{equation}
is exact, where $C^{^j}(\U,\P)=\underset{\alpha_0<\cdots<\alpha_j}
\prod\P(X_{\alpha_0}\cap\cdots \cap X_{\alpha_j})$ and $\delta$ is
given coordinatewise by the alternating sum of the restrictions. Notice that
$(\Cech{*}{\U}{\P},\delta)$ is a cohomological resolution of $\P(X)$.
\newline

For each perversity $\q$, the complex of $\q$-forms
$\Omintq{*}{-}$ is a presheaf on $X$ (and also on $B$).
The complex $\Inv{*}{-}$ of invariant
$\q$-forms is a presheaf on $X$ but
it is defined in the topology of invariant open sets;
we can also see it as a presheaf on $B$ up to a composition with
the orbit map. The complexes $\G{*}{-}$,
$\Resder{*}{-}$, and $\Resizq{*}{-}$
are presheaves on $B$.\newline

Because of the existence of controlled invariant partitions of the unity,
all these examples are sheaves.
We will study their cohomological properties on the charts
of $B$ and $X$.
\begin{lema}\label{lema de poincare}
    Let $\Phi:\s\times X\rightarrow X$ be a modelled action. Consider
    on $\R\times X$ the (obvious) modelled action trivial
    in $\R$. Then the projection $pr:\R\times X\rightarrow X$
    induces the following isomorphisms{\small
    \[
        \begin{array}{lclclcl}
        \Homintq{i}{\R\times X} & = & \Homintq{i}{X} & &
        \HG{i}{\R\times B} & =& \HGB{i} \\[0.2cm]
        \Hresder{i}{\R\times B}& = & \Hresder{i}{B} & &
        \Hresizq{i}{\R\times B} & =& \Hresizq{i}{B}
        \end{array}
    \]}
\end{lema}
\begin{proof} See \cite{gysin1}.
\end{proof}

\bprop\label{prop homgysin cono}
    Let $\Psi:\s\times L\rightarrow L$ be a modelled action on a
    compact stratified pseudomanifold $L$.
    For each perversity $\overline{0}\leq\q\leq\top$ and each $\epsilon>0$
    the map $\imath_\epsilon:L\rightarrow
    c(L)$ given by $p\mapsto [p,\epsilon]$
    induces the following isomorphisms
    {\small
    \begin{equation}
            \Homintq{i}{c(L/\s)}=\left\{
            \begin{array}{ll}
            \Homintq{i}{L/\s} & i\leq\q(\star) \\[0.2cm]
            0 & i>\q(\star)
            \end{array}
            \right.
    \end{equation}}
     Also
    {\small
    \begin{equation}
            \HG{i}{c(L/\s)}=\left\{
            \begin{array}{ll}
            \HG{i}{L/\s} & i\leq\q(\star)-2 \\[0.2cm]
            \ker\left[\partial:\HG{\q(\star)-1}{L/\s}
            \rightarrow \Homintq{\q(\star)+1}{L/\s}\right]
            & i=\q(\star)-1 \\[0.2cm]
            0 & i\geq\q(\star)
            \end{array}
            \right.
    \end{equation}}
    where $\partial$ is the connecting homomorphism of the Gysin sequence
    on  $L$.
\eprop
\bdem
    For the first isomorphism see \cite{illinois}.
    For the second one, we get a commutative diagram with exact horizontal rows
    {\small
    \[
        \begin{array}{ccccccccc}
            \rightarrow &                       
            \Homintq{i+1}{c(L/\s)} &
            \piarrow &
            \Homintq{i+1}{c(L)} &
            \rightarrow &
            \HG{i}{c(L/\s)}  &
            \bordarrow  &
            \Homintq{i+2}{c(L/\s)} &
            \rightarrow   \\
            & \downarrow &         
            & \downarrow &
            & \downarrow  &
            & \downarrow & \\
            \rightarrow &                                   
            \Homintq{i+1}{L/\s} &
            \piarrow &
            \Homintq{i+1}{L} &
            \rightarrow &
            \HG{i}{L/\s}&
            \bordarrow &
            \Homintq{i+2}{L/\s}  &
            \rightarrow
        \end{array}
    \]}
    where the vertical arrows are induced by $\imath_\epsilon:L\rightarrow c(L)$
    and $\imath_\epsilon:L/\s\rightarrow c(L/\s)$. For $i\leq\q(\star)-2$  we have
    enough vertical isomorphisms. By the Five Lemma
    \[
        \HG{i}{c(L/\s)}= \HG{i}{L/\s}
    \]
    For $i=\q(\star)-1$ the two left vertical arrows are isomorphisms.
    In the upper right coin we get $\Homintq{\q(\star)+1}{c(L/\s)}=0$.
    So
    \[
        \HG{\q(\star)-1}{c(L/\s)}=
        \text{\rm coker}(\pi^{*})=\ker\left[
        \HG{\q(\star)-1}{L/\s}
        \bordarrow \Homintq{\q(\star)+1}{L/\s}
        \right]
    \]
    For $i\geq\q(\star)$ the upper horizontal row has four zeros, thus
    $\HG{i}{c(L/\s)}=0$.
\edem
\bcor
    In the same situation of \S\ref{prop homgysin cono};
    if the vertex is not perverse then  the Gysin sequence of
    $c(L)$ is the Gysin sequence of $L$ truncated in
    dimension $i=\q(\star)-1$.
\ecor
\bprop\label{prop homresiduo cono}
    In the same situation of \S\ref{prop homgysin cono};
    if the vertex is a perverse stratum then
    the map $\imath_\epsilon:L\rightarrow
    c(L)$ induces the following isomorphisms
    {\small
    \begin{equation}\label{eq Hresdercono}
        \Hresder{i}{c(L/\s)}=\left\{
        \begin{array}{ll}
            \Hresder{i}{L/\s} & i\leq\q(\star)-3 \\[0.2cm]
            \ker\{\partial':
        \Hresder{\q(S)-2}{L/\s}\rightarrow\Hqm{\q(S)-1}{L/\s}\}
        & i=\q(\star)-2 \\[0.2cm]
            \ker
        \{\partial:\HG{\q(S)-1}{L/\s} \rightarrow \Homintq{\q(S)+1}{L/\s}\}
        & i=\q(\star)-1 \\[0.2cm]
            0 & i\geq\q(\star)
        \end{array}
    \right.
    \end{equation}}
    And{\small
    \begin{equation}\label{eq Hresizqcono}
        \Hresizq{i}{c(L/\s)}=\left\{
        \begin{array}{ll}
            \Hresizq{i}{L/\s} & i\leq\q(\star)-3 \\[0.2cm]
        \ker[\partial''\partial:\Hresizq{\q(\star)-2}{L/\s}
        \rightarrow \Homintq{\q(\star)+1}{L/\s}] & i=\q(\star)-2 \\[0.2cm]
            \displaystyle\frac{\Hqchi{\q(\star)-1}{L/\s}}{
        j^{*}(\ker^{^{\q(\star)-1}}(\partial))} & i=\q(\star)-1\\[0.2cm]
            0 & i\geq\q(\star)
        \end{array}
        \right.
    \end{equation}}
\eprop
\bdem
    We get the following commutative diagrams{\small
    \[
        \begin{array}{ccccccccc}
            \rightarrow &                       
            \Hqm{i}{c(L/\s)} &
            \overset{j^{*}}\rightarrow &
            \HG{i}{c(L/\s)} &
            \rightarrow &
            \Hresder{i}{c(L/\s)}  &
            \overset{\partial'}\rightarrow  &
            \Hqm{i+1}{c(L/\s)} &
            \rightarrow   \\
            (14) & \downarrow &         
            & \downarrow &
            & \downarrow  &
            & \downarrow & \\
            \rightarrow &                                   
            \Hqm{i}{L/\s} &
            \overset{j^{*}}\rightarrow &
            \HG{i}{L/\s} &
            \rightarrow &
            \Hresder{i}{L/\s}&
            \overset{\partial'}\rightarrow &
            \Hqm{i+1}{L/\s}  &
            \rightarrow
        \end{array}
    \]}
    and{\small
    \[
        \begin{array}{ccccccccc}
            \rightarrow &                       
            \HG{i}{c(L/\s)} &
            \overset{j^{*}}\rightarrow &
            \Hqchi{i}{c(L/\s)} &
            \rightarrow &
            \Hresizq{i}{c(L/\s)}  &
            \overset{\partial''}\rightarrow  &
            \HG{i+1}{c(L/\s)} &
            \rightarrow   \\
            (15) & \downarrow &         
            & \downarrow &
            & \downarrow  &
            & \downarrow & \\
            \rightarrow &                                   
            \HG{i}{L/\s} &
            \overset{j^{*}}\rightarrow &
            \Hqchi{i}{L/\s} &
            \rightarrow &
            \Hresizq{i}{L/\s}&
            \overset{\partial''}\rightarrow &
            \HG{i+1}{L/\s}  &
            \rightarrow
        \end{array}
    \]}
    where the horizontal rows are residual approximations and the
    vertical arrows are induced by the maps $\imath_\epsilon:L\rightarrow c(L)$
    and $\imath_\epsilon:L/\s\rightarrow c(L/\s)$.\newline

    Since the vertex is a perverse stratum $\pchi(\star)=1$ and
    $\m(\star)=2$. In the above diagrams the case $i\leq\q(\star)-3$ is direct from the
    Five Lemma and the case $i\geq\q(\star)$ is straightforward. We verify the
    cases $i=\q(\star)-2$ and $i=\q(\star)-1$ proceeding in two steps.\newline

    $\bullet$ {\tt Lower residue:} For $i=\q(\star)-2$, by \S\ref{prop homgysin cono}
    the two left vertical arrows in diagram (14) are isomorphisms.
    In the upper right corner we get
    $\Hqm{\q(\star)-1}{c(L/\s)}=0$. So {\small
    \[
        \Hresder{\q(\star)-2}{c(L/\s)} =
        \ker^{^{\q(\star)-2 }}
    [\partial':\Hresder{\q(\star)-2}{L/\s}
    \rightarrow\Hqm{\q(\star)-1}{L/\s}]
    \]}
    For $i=\q(\star)-1$ the upper corners are zeros. By
    \S\ref{prop homgysin cono} and the exactness of the upper
    horizontal row,
    \[
        \Hresder{\q(\star)-1}{c(L/\s)}=
    \HG{\q(\star)-1}{c(L/\s)}=
    \ker^{^{\q(\star)}}(\partial)
    \]

    $\bullet$ {\tt Upper residue:} For $i=\q(\star)-2$; by
    \S\ref{prop homgysin cono} the left vertical arrows in the diagram
    (15) are isomorphisms. Hence
    \[
        \imath_\epsilon:
    \Hresizq{\q(\star)-2}{c(L/\s)}\rightarrow
    \Hresizq{\q(\star)-2}{L/\s}
    \]
    is injective. We get a commutative exact diagram{\small
    \[
        \begin{array}{ccccccc}
        & 0  & & 0 & & 0 & \\
        & \downarrow & & \downarrow & & \downarrow & \\
        0 \rightarrow & \coker^{^{\q(\star)-2}}(j^{*},c(L/\s)) &
        \rightarrow & \Hresizq{\q(\star)-2}{c(L/\s)} &  \rightarrow &
        \ker^{^{\q(\star)-1}}(j^{\star},c(L/\s)) & \rightarrow 0 \\
        & \downarrow & & \downarrow^{\imath_\epsilon} & & \downarrow &  \\
        0 \rightarrow & \coker^{^{\q(\star)-2}}(j^{*},L/\s) &
        \rightarrow & \Hresizq{\q(\star)-2}{L/\s} &  \rightarrow &
        \ker^{^{\q(\star)-1}}(j^{\star},L/\s) & \rightarrow 0 \\
            & \downarrow & & \downarrow & & \downarrow & \\
        & 0 & \rightarrow & \coker^{^{\q(\star)-2}}(\imath_\epsilon) & \rightarrow &
        \frac{\ker^{^{\q(\star)-1}}(j^{\star},L/\s)}{
        \ker^{^{\q(\star)-1}}(j^{\star},c(L/\s))} & \rightarrow 0 \\
            & & & \downarrow & & \downarrow & \\
        & & & 0 & & 0 &
        \end{array}
    \]}
    So $\Hresizq{\q(\star)-2}{c(L/\s)}$
    is the kernel of the map{\small
    \[
        \overline{\imath_\epsilon}:
    \Hresizq{\q(\star)-2}{c(L/\s)}
    \rightarrow
        \frac{\ker^{^{\q(\star)-1}}(j^{\star},L/\s)}{
        \ker^{^{\q(\star)-1}}(j^{\star},c(L/\s))}
    =\frac{\Im^{^{\q(\star)-1}}(\partial',L/\s)}{
    \Im^{^{\q(\star)-1}}(\partial', L/\s)\cap
    \ker^{^{\q(\star)-1}}(\partial)}
    \]}
    In the last equality we used  \S\ref{prop homgysin cono},
    the exactness of the upper approximation and the fact that the
    third vertical arrow in  the diagram (15) is injective. So we can identify
    the image of the third vertical arrow with $\ker^{^{\q(\star)-1}}(\partial)$,
    the kernel of the connecting homomorphism of the Gysin sequence on $L$.
    We deduce that $\ker(\overline{\imath_\epsilon})$
    is the kernel of the composition
    \[
        \Hresizq{\q(\star)-2}{L/\s}
    \overset{\partial''}\longrightarrow
    \HG{\q(\star)-1}{L/\s}
    \overset{\partial}\longrightarrow
    \Homintq{\q(\star)+1}{L/\s}
    \]

    For $i=\q(\star)-1$ the first left vertical arrow in the
    diagram (15) is injective,
    the second is an isomorphism. In the upper right corner we get
    $\Hresizq{\q(\star)}{c(L/\s)}=0$.
    We obtain the exact commutative diagram {\small
    \[
        \begin{array}{ccccccc}
            0\rightarrow  & \coker^{^{\q(\star)-1}}(j,c(L/\s)) &
        \overset{\cong}\rightarrow & \Hresizq{\q(\star)-1}{c(L/\s)} &
        \rightarrow & 0 & \\
        & \downarrow^{\overline{\imath_\epsilon}} & &
        \downarrow^{\imath_\epsilon} & & \downarrow & \\
        0\rightarrow & \coker^{^{\q(\star)-1}}(j,L/\s) &
        \rightarrow & \Hresizq{\q(\star)-1}{L/\s} &
        \rightarrow & \ker^{^{\q(\star)}}(j,L/\s) &
        \rightarrow 0
    \end{array}
    \]}
    Notice that{\small
    \[
        \ker^{^{\q(\star)-1}}(\imath_\epsilon)
    \cong\ker(\overline{\imath_\epsilon})=
    \frac{\Im^{^{\q(\star)-1}}(j^{*},L/\s)}{
    \imath_\epsilon(\Im^{^{\q(\star)-1}}(j^{*},c(L/\s)))}
    =\frac{\Im^{^{\q(\star)-1}}(j^{*},L/\s)}{
    j^{*}(\Im^{^{\q(\star)-1}}(\imath_\epsilon))}
    =\frac{\Im^{^{\q(\star)-1}}(j^{*},L/\s)}{
    j^{*}(\ker^{^{\q(\star)-1}}(\partial))}
    \]}
    Also {\small
    \[
        \Im^{^{\q(\star)-1}}(\imath_\epsilon)=\coker^{^{\q(\star)-1}}(j^{*},L/\s)
    =\frac{\Hqchi{\q(\star)-1}{L/\s}}{
    \Im^{^{\q(\star)-1}}(j^{*},L/\s)}
    \]}
    So we get a short exact sequence{\small
    \[
        0\rightarrow
    \frac{\Im^{^{\q(\star)-1}}(j^{*},L/\s)}{
    j^{*}(\ker^{^{\q(\star)-1}}(\partial))}
        \rightarrow
    \Hresizq{\q(\star)-1}{c(L/\s)}
    \rightarrow
    \frac{\Hqchi{\q(\star)-1}{L/\s}}{
    \Im^{^{\q(\star)-1}}(j^{*},L/\s)}
        \rightarrow 0
    \]}
    We deduce that
    \[
        \Hresizq{\q(\star)-1}{c(L/\s)}=
    \frac{\Hqchi{\q(\star)-1}{L/\s}}{
    j^{*}(\ker^{^{\q(\star)-1}}(\partial))}
    \]
    This finishes the proof.
\edem
\bobss Statements \S7.2, \S7.3 and \S7.4 imply that
    \begin{itemize}
    \item[\Numero] For each fixed stratum $S$ in $X$ and
    each unfoldable chart
    $\beta:U\times c(L/\s)\rightarrow B$,
    ttatements we can calculate
    the cohomology of the Gysin term and the residues on
    the open $\Im(\beta)$.
    \item[\numero] According to \cite{borel},
    $\G{*}{-}$ is a constructible sheaf on $B$. Also
    $\Resder{*}{-}$ and $\Resizq{*}{-}$
    are constructible sheaves on $B$ with support on
    the perverse points' set $X^{perv}$.
    \end{itemize}
\eobss

\begin{teo}[The Gysin Theorem]\label{teorema gysin}
     Let $X$ be a stratified pseudomanifold, $\q$ a perversity in $X$,
     $\overline{0}\leq\q\leq\top$. For each modelled action
     $\s\times X\rightarrow X$ there are two
     long exact sequences relating the intersection cohomology of $X$ and $B$:
     the Gysin sequence
     \[
        \dots\rightarrow
        \Homintq{i+1}{X}\rightarrow
        \HGB{i}\bordarrow \Homintq{i+2}{B}
        \piarrow \Homintq{i+2}{X}
        \rightarrow\dots
     \]
     induced by the orbit map $\pi:X\rightarrow B$, and the lower approximation
     \[
    \dots\rightarrow
        \HGB{i}\rightarrow
        \Hresder{i}{B}\overset{\partial'}\rightarrow
        \Hqm{i+1}{B}\rightarrow
        \HGB{i+1}\rightarrow\dots
     \]
     induced by the inclusion $\Omint{\q-\m}{*}{B}\iarrow \GB{*}$. These sequences
     satisfy
     \begin{itemize}
     \item[\Numero] If $X$ has no perverse strata then
     $\GB{*}=\Oqchi{*}{B}=\Oqm{*}{B}$, $\Resder{*}{B}=0$ and the
     connecting homomorphism of the Gysin sequence is the multiplication
     by the Euler Class $\varepsilon\in \Homint{\pchi}{2}{B}$.
     \item[\numero] If $X$ has perverse strata, then
     $\Hresder{*}{B}$ is calculated through a spectral sequence
     in $B$, whose second term
     \[
         E_{_2}^{^{ij}}=\Hcech{j}{X^{perv}}{\P^{^{i}}}
     \]
    is the cohomology of the perverse points' set $X^{perv}$
    with values on a locally constant graduated constructible presheaf $\P^{^{*}}$.
    For each fixed point $x\in X$ the stalks
    {\small
    \[
        \P^{^i}_{_x}=\left\{
        \begin{array}{ll}
            \Hresder{i}{L/\s} & i\leq\q(S)-3 \\[0.2cm]
            \ker\{\partial':
        \Hresder{\q(S)-2}{L/\s}\rightarrow\Hqm{\q(S)-1}{L/\s}\}
        & i=\q(S)-2 \\[0.2cm]
            \ker
        \{\HG{\q(S)-1}{L/\s} \overset{\partial}\rightarrow \Homintq{\q(S)+1}{L/\s}\}
        & i=\q(S)-1 \\[0.2cm]
            0 & i\geq\q(S)
        \end{array}
    \right.
    \]}
    depend on the Gysin sequence and the residual approximation
    induced by the action of $\s$ on of the link $L$
    of the stratum containing $x$.
    \end{itemize}
\end{teo}
\bdem
    Statement (1) has been already  proved in the preceding sections.
    Statement (2) arises from the usual spectral sequence induced by a double
    complex; see for instance \cite{bott}, \cite{godement}.
    The double complex we take is the residual Cech double complex
    \[
        (\Cech{j}{\U}{\Resder{i}{-}},\delta,d)
    \]
    induced by an invariant open cover $\U=\{B_\alpha\}_\alpha$ of $B$;
    where $\delta$ is the Chech differential induced by the restrictions, and
    $d$ is the usual differential operator. We define
    $\U$ as follows: First take the complement of the fixed
    points' set $X_0=X-X^{\s}$; we ask $B_0=\pi(X_0)$ to be in $\U$.
    Second, for each fixed point
    $x\in X$ we take a modelled chart
    \[
        \alpha:U_\alpha\times c(L)\rightarrow X
    \]
    in the stratum $S$ containing $x$; such that $x\in U_\alpha$. We ask the
    $U_\alpha$'s intersecting $S$ to be a good cover of $S$. We take
    \[
        B_\alpha=\pi(\Im(\alpha))\cong U_\alpha\times c(L/\s)
    \]
    Since the sheaf $\Resder{*}{-}$ vanishes identically on $B_0$; the second term
    of the spectral sequence is
    \[
        E_{_2}^{^{ij}}=
    H^{^j}_{_\delta}H^{^i}_{_d}(\U,\Resder{}{-})=
    \Hcech{j}{\U}{\mathcal{H}\Resder{i}{-}}=
    \Hcech{j}{X^{perv}}{\mathcal{H}\Resder{i}{-}}
    \]
    So $\P^{^{*}}=\mathcal{H}\Resder{*}{-}$ is the desired presheaf.
    The remarks on the stalks are immediate from
    \S\ref{lema de poincare}, \S\ref{prop homgysin cono}.
\edem

\bfalsa{\bf Exceptional actions}
A modelled action $\Phi:\s\times X\rightarrow X$ is {\bf exceptional}
if the links of $X$ have no perverse strata; i.e., if any perverse stratum
of $X$ is a closed (minimal) stratum.
\efalsa
\bcor
    For any exceptional action $\Phi:\s\times X\rightarrow X$
    we have
    \begin{equation}\label{eq resdermin}
        \Hresder{*}{B}=
    \underset{S}\prod\ \Hcech{*}{S}{\mathfrak{Im}_{_{\q}}(\varepsilon_L)}
    \end{equation}
    where $S$ runs over the perverse strata and
    $\Hcech{*}{S}{\mathfrak{Im}_{_{\q}}(\varepsilon_L)}$ is the cohomology of $S$
    with values on a locally constant presheaf
    with  stalk
    {\small\[
         \ker\{\varepsilon_L:
     \Hqchi{\q(S)-1}{L/\s}\rightarrow
     \Homintq{\q(S)+1}{L/\s}\}
    \]}
    the kernel of the multiplication by the Euler class $\varepsilon_L\in
    \Homint{\pchi}{2}{L/\s}$ of the link $L$ of $S$.
\ecor
\bdem
    If the link $L$ of a perverse stratum $S$ has no perverse strata, then
    the stalk $\mathcal{H}^{^{i}}\Resder{}{-}$ vanishes for
    $i\neq \q(S)-1$ (so it is a single presheaf). The equality (\ref{eq resdermin})
    is straightforward, since the perverse strata are disjoint closed subsets.
\edem

\section*{Acknowledgments}
We would like to thank the accurate observations of
M. Saralegi. While writing this article, the author received
the financial support of the CDCH-Universidad Central de
Venezuela and the hospitality of the staff in the
Math Department-Universit\'e D'Artois.

\end{document}